\documentclass[12pt]{amsart}
\usepackage{amsfonts,amsmath,amsthm,amscd,amssymb,latexsym,cite,verbatim,texdraw,floatflt,
caption2}
\RequirePackage{hyperref}

\usepackage[a4paper]{geometry}

\theoremstyle
{plain}

\begin{document}

\title{The dynamical approach to the conjugacy in groups      }

\author{ Igor Protasov, Ksenia Protasova}

\maketitle
\vskip 5pt

{\bf Abstract.}  
Given a discrete group $G$, we identify the Stone-$\check C$ech compactification  $\beta G$ with the set of all ultrafilters on $G$ and put $G^\ast =\beta G\setminus G$. The action $G$ on $G$ by the conjugations $(g,x)\mapsto g^{-1}xg$ induces the action of $G$  on $G^\ast$  by $(g, p)\mapsto  p^g $, $p^g = \{ g^{-1} Pg: P\in p\}$. We study interplays between the algebraic properties of $G$ and the dynamical properties of $(G, G^\ast)$. In particular, we show that $p^G$ is finite for each $p\in G^\ast$ if and  only if the commutant of $G$ is finite.

\vspace{6 mm}

 1991 MSC: 20E45, 54D80.

\vspace{3 mm}

Keywords:  space $G^\ast$ of ultrafilters on a  group $G$,  the action of $G$ on $G^\ast$ by the conjugations, ultracenter of $G^\ast$.

\section{Introduction and preliminaries  }

Let a group $G$ acts on a set $X$ by 
$$G\times X \longrightarrow X: (g,x)\mapsto gx. $$

We endow $X$ with the discrete topology, identity the Stone-$\check{C}$ech extension $\beta X$ of 
$X$ with the set of all ultrafilters on $X$  and put $X^\ast = \beta X\setminus X$. 
By the universal property of $\beta X$, every mapping $f: X\longrightarrow K$, where $K$ is a compact Hausdorff space, gives up to the continuous mapping $f^\beta : \beta X\longrightarrow K $. 
Using this property, we extend the action of $G$  on $X$  to the action
$$G \times \beta G \longrightarrow \beta G : (g,p)\mapsto gp ,$$ 
where $gp$ is the ultrafilter $\{ gP: P\in p \}$. We note that $gp\in X^\ast$ for each $g\in G$, $p\in X^\ast$.
Thus, we get the {\it dynamical system}   $(G, X^\ast).$

In the case $X=G$ and $gx$ is the product of $g$ and 
$x$
in $G$, $(G, G^\ast)$ is studied in \cite{b2}.
The general case of $(G, X^\ast)$ aroused in {\it Asymptology}, see the survey \cite{b4} and [1, Section 5].

In this paper, $X=G$, but $G$ acts on $X$ by the conjugations 
$(g,x)\mapsto g^{x}$, 
$g^{x}= x^{-1} gx$.

In what follows, all groups under consideration are suppose  to be infinite. If $g\in G$ and $p\in G^\ast$ then $$ p^g =\{ g^{-1} Pg: P\in p\}, \ \  p^G=\{ p^g: g\in G\}$$ and $p^G$  is called the $G$-{\it orbit} of $p$.

We note that the closure $cl p^G$ of $p^G$
coincides with the set $\{p^q : q\in \beta G \}$,
where 
$p^q =q$-$lim \{ p^x : x\in G\}$.
To describe $p^q$  more explicitly, we take $Q\in q$ and, for each $x\in Q$, pick $P_x \in p$. 
Then  $\bigcup \{ x^{-1} P_x x : x\in Q \} \in p^q$ 
and these subsets from a base of $p^q$.

For a group $G$ and $p\in G^\ast$,
 we denote 
$$ St (p) = \{ g\in G: p^g =p  \}, \ \  Z_G = \{ g\in G: g^x =g \  \text{ for \ each }\ x\in G\} ,$$

$$UZ_G = \{ p\in G^\ast : St(p)=G\}.  $$ 
and say that $UZ_G$ is the {\it ultracenter} of $G^\ast$. 

We note that $g\in St(p)$ if and  only if there exists $P\in p$ such that $gx=xp$ for each $x\in P$, and $Z_G^\ast\subseteq UZ_G$.
We use also the following observation:  $|p^G|= |G: St(p)|$.

\vskip 10pt

{\bf Proposition 1.1.} {\it 
A group $G$ is Abelian if and only if 
$UZ_G = G^\ast$.

\vskip 7pt

Proof.  } 
If $G$ is Abelian then $St(p)=G$ for each $p\in G^\ast$ so  
$G^\ast \subseteq   UZ_G $.

We suppose that $G^\ast = UZ_G$.
If $g$, $x\in G$, $p\in G^\ast$ then $(xp)^g=xp$ and 
$(xp)^g=x^g p^g =x^g  p$ 
so $x^g =x$ and $G$ is Abelian. $ \ \ \Box$ 

\vspace{10 mm}

The multiplication on a discrete group $G$ can be extended (see \cite{b2}) 
on $\beta G$  so that, for  each  $g\in G$,
the mapping

$$\beta G \longrightarrow \beta G: p\mapsto gp $$
is continuous, and for each $p\in \beta G$, the  mapping 
$$\beta G \longrightarrow \beta G: q\mapsto qp $$
is continuous. Then,
$\beta G$ is a right topological  semigroup and $G^\ast$  is a subsemigroup of $\beta G$.

For $q, p \in \beta G$, we choose $Q\in q$ and, for each $g\in Q$, pick $P_g\in p$. Then $\cup \{x P_x : x\in Q \} \in qp$  and subsets of this form are a base for $q,p$. 
   
\vskip 10pt

{\bf Proposition 1.2.} {\it 
For a group $G$,
either $UZ_G$ is empty or $UZ_G$  is a subsemigroup of $G^\ast$.
\vskip 7pt

Proof.  } 
It suffices to note that $(qp)^g = q^g p^g$ for all $g\in G $ and $q, p\in G^\ast$. $ \ \ \Box$ 
 
In Section 2, we show that, for every  group $G, G^\ast$ has a plenty of discrete orbits $p^G$ and if $G$ is countable then the set $\{ \ p\in G^\ast : p^G$ is discrete $\}$ is dense in $G^\ast$. We give some sufficient conditions on a subset $A$ of $G$ under which some orbit $p^G$, $p\in A^\ast$ or each orbit $p^G, \ p\in A^\ast$ is infinite. 

In Section 3, we characterize groups $G$ such that each orbit $p^G$,  $p^{G^\ast}$ is finite.

In Section 4, we show that every dynamical system $(H, X^\ast)$ can be embedded into $(G, G^\ast)$  for an appropriate choice of $G$. 
This means that the dynamical systems of the form $(G, G^\ast)$  could be  very complicated. 


\section{ Plenty of orbits}
    

{\bf Theorem 2.1.} {\it 
For every subset $A$ of a group $G$ such that 
$|A|=|G|$, we have }
$$ |\{ p\in G^\ast : A\in p \text{  and  } p^G \text{ is discrete } \}|= 2^{2^{|A|} } .$$

\vskip 7pt

{\it Proof.  } We enumerate $G=\{ g_\alpha : \alpha < |G|\}$, put $K_\alpha = \{g_\beta : \beta<\alpha \}$  and choose inductively $\{a_\alpha : \alpha < |G| \}$ such that 
$$ a_\alpha ^{K_\alpha} \bigcap  a_\beta ^{K_\beta} =\emptyset $$  for all $\alpha < \beta$.
If $g\in G$, $p\in G^\ast$, $\{ a_\alpha : \alpha < |G|  \}\in p$ and $|P|=|G|$ for each $P\in p$  then either $p^g=p$ or  
$( G\setminus \{ a_\alpha : \alpha < |G|  \}) \in p^G$.
It suffices to note that we have $2^{2^{|G|}}$ possibilities to choose $p$ with these properties. $ \ \ \Box$

\vskip 10pt

{\bf Theorem 2.2.} {\it 
For  a countable group $G$, the following 
statements hold 
 
 \vskip 7pt

$(i)$  the set   $\{ p\in G : p^G \ is \ discrete \}$ 
 is dense in 
 $G^\ast$;
\vskip 7pt

 $(ii)$  if  $cl \ p^G\bigcap cl \ q^G \neq \emptyset $   then either 
$cl \ p^G\subseteq cl \ q^G $ or $cl \ q^G\subseteq cl \ p^G $;
\vskip 7pt

 $(iii)$  there exist an infinite subset $A\subseteq G$ such that 
 $cl \ p^G\bigcap cl \ q^G = \emptyset $  for all distinct $p, q \in A^\ast$.
 
 }

\vskip 10pt

{\it Proof.  } 
$(i)$  Apply Theorem 2.1.

\vskip 7pt

$(ii)$ Apply the Frolik's lemma [2, Corollary 3.42];
\vskip 7pt

$(iii)$ We take  a subset $A$ from the proof of  Theorem 2.1. Let $p, q \in A^\ast$, 
$p\neq q$
 but 
 $cl \ p^G\bigcap cl \ q^G = \emptyset $. 
By (ii), either $p=cl \ q^G$ or $q=cl \ p^G$.
In both cases, there exists
 $r\in A^\ast $,  
 $g\in G$
such that 
$r^g \neq r$ 
and 
$r^g \in A^\ast$, 
contradicting the choice of $A$. 
$ \ \ \Box$     

\vskip 10pt

{\bf Theorem 2.3.} {\it If some conjugate class $a^G$ of a group $G$ is infinite then there exists $p\in G^\ast$ such that $p^G$ is infinite.

\vskip 7pt

Proof. } We choose a sequence $(g_n)_{n\in \omega}$ in $G$ such that $a^{g_n}\neq a^{g_m}$ for all distinct $m, n$.
Let $q\in G^\ast$. If there exists an infinite subset $I\subseteq \omega$ such that $q^{g_m} \neq q^{g_n}$  for all distinct $m, n$
then we put $p=q$.

Otherwise, there exists an infinite subset $J\subseteq \omega $ such that $q^{g_n} = q^{g_m}$
for  all distinct $m, n \in J$. 
Then $(aq) ^{g_n} = a^{g_n} q$, $(aq) ^{g_m} = a^{g_m} q$.
Since $a ^{g_n} \neq a^{g_m} $, we have 
$(aq) ^{g_n} \neq (aq) ^{g_m}$, put $p=aq$.
$ \ \ \Box$

\vskip 10pt

For a subset $A$ of a group $G$ and $p\in G^\ast$, we denote $\Delta_p (A)= p^G \bigcap A^\ast$.

\vskip 10pt

{\bf Theorem 2.4.} {\it 
Let $A$ be a subset of a group $G$ such that 
$\Delta_p (A)$ is infinite  for each $p\in A^\ast$. 
Then there exists $q\in A^\ast$ such that $q^G$ is not discrete. 

\vskip 7pt

Proof. } We use the Zorn's lemma and compactness of 
$G^\ast $ to choose the minimal by inclusion, closed subset $S$ of $G^\ast $ such that $S\cap Y^\ast\neq \emptyset$ and  $g^{-1} S g \subseteq S$ for each  $g\in G$. Let $q\in S \cap Y^\ast$. 
Since 
$q^G \cap Y^\ast$ is infinite, there exists a limit point 
$r\in Y$  of  $q^G $.
By the minimality of $S, q$ is a limit point of $r^G$.
It follows  that $q$ is not isolated in $q^G $ and $q^G $ is not discrete. $ \ \ \Box$ 

\vskip 10pt

Let $(y_n)_{n\in\omega}$,  $(b_n)_{n\in\omega}$ be  sequences in a group $G$.
Adopting the dynamical terminology from  \cite{b5}  to the action of $G$ on $G$ by conjugations, we say that 
\begin{Large}
  $$ Y=  \Big\{ b_n ^{y_n ^{\varepsilon_n}\dots y_0 ^{\varepsilon_0}} : n<\omega, \  \varepsilon_i \in \{0,1\} \Big\} $$
\end{Large}
is a {\it piece-wise shifted FP-set}  if

\begin{Large}
  $$ (*) \  b_n ^{y_n ^{\varepsilon_n}\dots y_0 ^{\varepsilon_0}} \ = \ b_m  ^{y_m ^{\delta_n}\dots y_0 ^{\delta_0}} ,   \  \   \   \varepsilon_i , \delta_j\in  \{0,1\} $$ 
\end{Large}
implies $n=m $ and $\varepsilon_0 = \delta_0 , \dots , \varepsilon_n = \delta_n$.

\vskip 10pt

{\bf Theorem 2.5.} {\it If $Y$ is a piece-wise shifted FP-set in a group $G$ and 
$p\in Y^\ast$ then the orbit $p^G$ is infinite.
\vskip 7pt

Proof. } 
For $n<\omega$, we put 

\begin{Large}

  $$ Y_m  =  \Big\{ b_n ^{y_n ^{\varepsilon_n}\dots y_0 ^{\varepsilon_0}} : \varepsilon_0 =  \dots = \varepsilon_{m-1} =0, 
 \ \  \varepsilon_m =1  \Big\} $$
\end{Large}
and note that $Y_m \cap Y_m = \emptyset$ if  $m\neq m^\prime$.

If $Y_m \notin p$ for each $m<\omega $ then  
$p^{y_n} \neq p^{y_k} $ for  all distinct  $n, k$.

If $Y_m \in p$  then $p^g \neq p $ for $g=y_m$ and we repeat  the arguments for $p^g$ in place of $p$. $ \ \ \Box$ 

\vskip 10pt

\section{Finite orbits }

We recall that $G$ is an FC-group if $a^G$ is finite for $a\in G$.

\vskip 10pt

{\bf Theorem 3.1.} {\it For a group  $G$, the following conditions are equivalent

\vskip 7pt

{(i)}  $p^G$ is finite for each $p\in G^\ast$;
\vskip 7pt

 {(ii)}  there exists a natural number $n$ such that $|p^G|< n$ for each $p\in G^\ast$;
\vskip 7pt

 {(iii)}  there exists a natural number $n$ such that $|a^G|< n$ for each $a\in G$;
 
 \vskip 7pt
 
 {(iv)}   the commutant   of $G$  is finite.
 
\vskip 10pt

Proof. } 
The equivalence $(iii)$ $\Leftrightarrow  $ $(iv)$ is proved in \cite{b3},
$(ii)$ $\Rightarrow  $ $(i)$ is evident.
\vskip 7pt

$(i)$ $\Rightarrow  $ $(iii)$.
If some class $a^G$ is infinite then, by Theorem 2.3, some orbit $p^G$ is infinite. 
We suppose that $G$ is an FC-group, but for every natural number $n$, there exists $a\in G$ such that $|a^G|>n$. 
In view of Theorem 2.5, it suffices to construct a piece-wise  shifted FP-set $Y$ in $G$.

We choose $c_0 , y_0 \in G$ such that  $c_0^{y_0} \neq c_0$, put $b_0 = c_0$.
Suppose that we have chosen $y_0, \dots , y_k$   and 
$c_0, \dots , c_k$ such that, for $b_i = c_0, \dots , c_i$, $(\ast)$ holds  for all $n, m \leq k$. 
We denote by $Z$ the intersection 
of the centralizers of
$c_0 , \dots , c_k$, $y_0 , \dots , y_k$.
Since $G$ is an FC, $|G: Z|$ is finite. 
It follows that $Z$ has an arbitrary large class of conjugate  elements in $Z$. Using this observation, we can choose $c_{k+1}, y_{k+1}\in Z$  such that, for $b_{k+1}= b_k c_{k+1}$, $(\ast)$ holds for all $n,m \leq k+1$.

After $\omega$ steps,  we  get the desired piece-wise shifted FC-group
\begin{Large}
$$Y= \{ b_n ^{y_n ^{\varepsilon_n} , \dots , y_0 ^{\varepsilon_0} } : n<\omega, \varepsilon_i \in \{0,1\} \}. $$
    
\vskip 7pt

\end{Large}

$(iii)$ $\Rightarrow  $ $(ii)$.
We suppose that there exists $p\in G^\ast$ such that $|p^G|\geq m$.
We choose $g_1, \dots , g_m$
in $G$ such $p^{g_i} \neq p^{g_j} $ 
for all $i<j\leq m$.
Then there exists $P\in p$ such that 
$P^{g_i} \bigcap P^{g_j}=\emptyset$
for all $i<j\leq m$.
If $a\in P$ then $a^{g_i} \neq a^{g_j} $ for 
$i<j\leq m$, contradicting $(iii)$.
$ \ \ \Box$

\vskip 10pt

{\bf Example 3.2.} We define a transitive group $S$ of permutations of $\omega$ such that each orbit $Sp$, $p\in \omega^\ast$ is finite, 
and for every natural number $n$ there exists $q\in\omega^\ast$ such that $|Sq|=n$.  

 First, partition  $\omega$  infinite subsets 
$\{ W_n : n\in \omega\}$.
Second, partition each $W_n$ into infinite subsets $\{ W_{ni} : i< n \}$.  Third, for each  $n$, defines a bijection 
$b_n$ of $\omega$  such that 
$f_n (W_{n0})= W_{n1}, \ $  
$ \ f_n (W_{n1})= W_{n2},\dots ,\ $   $ \ f_n (W_{n \ {n-1}})= W_0, \ $ 
and $ \ f_n (x)=x$ for each $x\in \omega \setminus W_n$.
At last, take a group  of permutations of $\omega$
generated by $\{f_n : n<\omega \}$ and all permutations with finite support. 
It $p\in \omega\ast$ and $W_n \in p$ then $|Sp|=n$, 
otherwise, $|Sp|=1$.

\vskip 10pt

{\bf Question 3.3.} {\it Let   $G$ be a group such that 
 $p^G$
is discrete for each $p\in G^\ast$. 
Is $p^G$ finite for each $p\in G^\ast$?}

\vskip 7pt

In light of Theorems 2.4 and 2.5, to answer this question in affirmative, it suffices to show that if some orbit $p^G$ is infinite then $G$ contains a piece-wise shifted FP-set.

\vskip 10pt

\section{Dynamical embeddings }

Let $(H, Y^\ast)$, $(G, X^\ast)$
be dynamical systems. We suppose that there exists an isomorphic embedding  $f: H \longrightarrow G$ and an injective mapping $\phi: Y\longrightarrow X$ 
such that, for each $(h, p)\in (H, Y^\ast)$, we have 
$\phi^\beta (hp)= f(h) \phi^\beta (p))$.
We say that the pair $(f,\phi)$ is a dynamical embedding of $(H, Y^\ast)$  into $(G, X^\ast)$.

\vskip 10pt

{\bf Theorem 4.1.} {\it Every dynamical system $(H, Y)$ admits a dynamical embeddings into $(G,G^\ast)$ 
for an appropriate choice  of a group $G$.

\vskip 7pt

Proof.} We consider $\{0, 1 \}^Y$
as a group with 
 point-wise 
addition $mod \ 2$. 
For $h\in H$ and $\chi \in \{0, 1 \}^Y$, 
$\chi_h$
is defined by  $\chi_h (y)= \chi (h^{-1}y)$.
Then we define a semidirect product 
$G =  \{0, 1 \}^Y  \leftthreetimes  H$ by

$$(\chi, h) (\chi^\prime, h^\prime) = (\chi + \chi_h 
^\prime , \  hh^\prime).$$

Let $f: H\longrightarrow G$
be the natural 
 isomorphic embedding and $\phi : Y \longrightarrow G$ is defined by 
 $\phi (y) = ( \psi _y, id)$, $\psi _y$
is the characteristic function of $\{y\}$, $id$ is the identity of $H.$ 
Clearly, 
$\phi ^\beta (hp)= (\phi ^\beta (p))^{f(h)}  $
so $(f,\phi)$ is the desired dynamical embedding. 
$ \ \ \Box$

\vskip 10pt

{\bf Example 4.2.} Let $S$ be the group of all permutations of $\omega = \{ 0,1, \dots \} $.
Every orbit $Sp$, $p\in \omega^\ast$ is not discrete. 
Applying  Theorem 4.2, we conclude that Theorem 2.2(i) needs not to be true  for an uncountable  group.

\vskip 10pt

{\bf Example 4.3.}
Let $F$ be the group of all permutations of $\omega$ with finite supports. Then $Fp= \{p \}$  for  each $p\in \omega^\ast$. We use $G$ given by Theorem 4.1 for 
$(F, \omega^\ast)$.
Then $ \phi^\beta  (\omega^\ast) \subseteq U Z _G$.


\vskip 10pt

CONTACT INFORMATION
\vskip 10pt

I.~Protasov: \\
Faculty of Computer Science and Cybernetics  \\
        Kyiv University by Taras Shevchenko \\
         Academic Glushkov pr. 4d  \\
         03680 Kyiv, Ukraine \\ i.v.protasov@gmail.com

\vskip 10pt

K.~Protasova: \\
Faculty of Computer Science and Cybernetics  \\
        Kyiv University \\
         Academic Glushkov pr. 4d  \\
         03680 Kyiv, Ukraine \\ k.d.ushakova@gmail.com

\end{document}